\documentclass[12pt,a4paper]{article}
\usepackage{amssymb,amsmath,amsthm}
\usepackage{graphicx, color}
\def\de{{\rm d}}
\def\dA{\delta_{\theta_1}}
\def\dB{\delta_{\theta_2}}
\def\dAA{\delta_{\theta_1}^2}
\def\dBB{\delta_{\theta_2}^2}

\def\E{{\rm E}}

\newtheorem{theorem}{Theorem}[section]

\title{Parametric estimation for partially hidden diffusion processes  
sampled at discrete times}
\author{Iacus S.M.\footnote{Department of Economics, Business and  
Statistics, University of Milan, Via Conservatorio 7, 20122 Milan,  
Italy},
Uchida M.\footnote{{Graduate School of Engineering Science, Osaka  
University, Toyonaka, Osaka 560-8531, Japan }},
Yoshida N.\footnote{Graduate School of Mathematical Sciences,  
University of Tokyo, 3-8-1 Komaba, Meguro-ku, Tokyo 153-8914 Japan}}
\begin{document}

\date{}
\maketitle

\begin{abstract}
For a one dimensional diffusion process $X=\{X(t) \ ; \ 0\leq t \leq T \}$,
we suppose that $X(t)$ is hidden if it is below some fixed and known  
threshold $\tau$, but otherwise it is visible.
This means a partially hidden diffusion process.
The problem treated in this paper
is to estimate finite dimensional parameter in both drift and  
diffusion coefficients under
a partially hidden diffusion process obtained by a discrete sampling  
scheme.
It is assumed that the sampling occurs at regularly spaced time  
intervals of length $h_n$
such that $n h_n=T$.
The asymptotic is when $h_n\to0$, $T\to\infty$ and $n h_n^2\to 0$ as  
$n\to\infty$.
Consistency and asymptotic normality for estimators of parameters in  
both drift and diffusion coefficients {are} proved.
\end{abstract}



\newpage

\section{Introduction}
We consider the estimation of the unknown parameter $\theta= 
(\theta_1, \theta_2)$ characterizing {{a}}
{one-dimensional diffusion process defined by
the stochastic differential equation}
$$
dX(t) = b(X(t), \theta_2) d t + \sigma(X(t), \theta_1) dW_t, \quad X 
(0)=x_0, \quad t \in [0,T],
$$
where $W$ is a {one-dimensional} standard Brownian motion,
$b$ and $\sigma$ are supposed to be regular enough to ensure the  
existence of a {(strong)} solution to the above stochastic  
differential equation.
In the situation where discrete observations are
${\bf X}_n =\{ X({t_i}) \ ; \ i=0,1,\ldots, n \}$ with $t_i = i h_n$,  
$n h_n = T$, the estimation problem for\\ the parameter $\theta$
has been considered by several authors,
see Florens-Zmirou (1989), Prakasa-Rao (1983, 1988), Yoshida (1992,  
2005) and Kessler (1997).
In this paper, however, we generalize it to a different setup.
We suppose that $X(t)$ is observable if $X(t) > \tau$ for some  
threshold $\tau$,
and that $X(t)$ can not be observed if $X(t) \leq \tau$.
This means that the original process becomes a partially hidden  
diffusion process based on a threshold $\tau$,
and the discretized trajectory ${\bf X}_n$
is also influenced by a threshold $\tau$. 
This type of observation naturally  arises in the study of stochastic  
resonance and has been treated
so far in the statistical context for the i.i.d. case in Greenwood  
{\it et al.} (2000), for continuous time
{ergodic} diffusion processes in Iacus (2002) and for a class of  
continuous time mixing processes in Iacus and Negri (2003).
In signal theory this corresponds to the problem of signal detection  
when the signal is so faint that it is not always receivable by some  
detector. This scheme of observation  frequently appears in radio and  
CCD  astronomy in the problem of identification of faint perturbed  
signals originated by astronomical sources (see e.g. Starck {\it et  
al.}, 1999). Partially observed diffusion model also arises in the  
context of financial markets (see e.g. Zeng, 2003) and in neuronal  
activation analysis (see e.g. Movellan and  Mineiro, 2002). In  
stochastic resonance context the original observation is altered by  
adding some noise with known structure to the channel in order to  
have full (but eventually quite noisy) observations, hence the  
problem is the one of determining the optimal level of noise. In the  
approach used in this paper, only the available observations are  
retained and used to estimate $\theta$.
In this setup, we need to build a  contrast function which is  
different from the one proposed in the literature of estimation for  
discretely observed diffusion processes cited above.
Other different {approaches} based on particle filters (see e.g.  
Fearnhead {\it et al.}, 2006)
and observation augmentation (see e.g. Roberts and Stramer, 2001)
{have} been also recently proposed in the literature
but our approach and asymptotic scheme adopted is substantially  
different from these references.
Nevertheless, after some refinement it is still possible to prove  
consistency and asymptotic normality of the proposed estimators
along the lines of e.g. Yoshida (1990, 1992),
Genon-Catalot and Jacod (1993)
and Kessler (1997).
The organization of the paper is as follows: Section \ref{sec:model}  
introduces the model, the assumptions
and two contrast functions.
of observation.
Section \ref{sec:main} contains the statement of the main result on  
consistency and asymptotic normality of estimators.
Section \ref{sec:proofs} is devoted to the proofs of the results in  
Section \ref{sec:main}.


\section{Model of observation and assumptions}\label{sec:model}
{Let}
$X = \{X(t) \ ; \ 0 \leq t\leq T\}$
{denote}
{{a}} diffusion process
{{satisfying}}
\begin{equation}
d X(t) = b(X(t), \theta_2) d t + \sigma(X(t), \theta_1) d W_t, \quad X 
(0) = x_0, \quad t \in [0,T].
\label{eq:sde}
\end{equation}
The parameter of our interest is $\theta= (\theta_1, \theta_2)$, $ 
\theta\in\Theta$
and $\Theta$ is a compact {{rectangle}} in ${\bf R}^2$.
The true value is denoted by $\theta_0= (\theta_{1,0}, \theta_{2,0})$
and it is assumed that $\theta_0 \in \mbox{Int}(\Theta)$.
Let $X_i = X(t_i)$,
$t_i = i h_n$ and $n h_n = T$.
For $i=0,1,\ldots, n$, we assume that $X_i$ is observable if $X_i >  
\tau$ for some threshold $\tau$,
and that $X_i$ is unobserved if $X_i \leq \tau$.
The asymptotics 
will be investigated when $h_n\to 0$, $nh_n \rightarrow \infty$ and  
$n h_n^2\to 0$ as $n \to \infty$.
In order to simplify the description, we use the following notation
$$
\sigma_{i} = \sigma(X_i,\theta_1), \quad b_{i} = b(X_i,\theta_2),  
\quad \Delta_i X = X_i - X_{i-1}.
$$
When the coefficients are evaluated at the true value of the  
parameter, we will write
$$
\sigma_{i}^* = \sigma(X_i,\theta_{1,0}), \qquad b_{i}^* = b(X_i, 
\theta_{2,0}).
$$
We further
define
$\delta_{\theta_i} f= \frac{\partial}{\partial\theta_i} f$.
For any real sequence {$u_n \in (0{{,}}1]$},
$R(u_n,x)$ represents a function such that
\begin{equation}
{|R(u_n, x)|} \leq u_n C(1+|x|)^C,
\label{eq:Rn}
\end{equation}
where $C$ is a {positive} constant independent of $n$ and $x$ (and  
eventually $\theta$ when $x$ is $X(t)$).
In the proof, $K$ and/or $C$ denote generic constants
independent of $\theta$, $x$ and $n$.

\subsection*{Assumptions}

\begin{enumerate}
\item there exists {$K>0$} such that
$$|b(x,\theta_{2,0}) - b(y,\theta_{2,0})| +  |\sigma(x,\theta_{1,0})  
- \sigma(y,\theta_{1,0})|\leq K|x-y|,$$
so that  \eqref{eq:sde} has a unique solution for $\theta=\theta_0$.  
\label{hip1}
\item
{the process $X$ is stationary and ergodic for $\theta=\theta_0$
with its invariant measure denoted by $\nu_{\theta_0}$}. \label{hip2}
\item
{for all $p\geq 0$,
$\sup_t \E |X(t)|^p<\infty$}. \label{hip3}
\item $\inf_{x,\theta_1} \sigma^2(x,\theta_1) = K_4 >0${.} \label{hip4}
\item (polynomial growth) the coefficients $b$ and $\sigma$ are  
{continuously} differentiable with respect to $x$ up to order 2
for all $\theta_1$ and $\theta_2$.
Themselves and their derivatives up to order 2 are of polynomial  
growth in $x$, uniformly in $\theta$. \label{hip5}
\item (polynomial growth) the coefficients $b$ and $\sigma$ and all  
their $x$ derivatives up to order 2,
are three times {continuously} differentiable with respect to $\theta 
$ for all $x$. Moreover, these $\theta$-derivatives are of polynomial  
growth in $x$ and uniformly on $\theta$. \label{hip6}


\item (identifiability)
$b({{\cdot}},\theta_2) = b({{\cdot}},\theta_{2,0})$
and
$\sigma({{\cdot}},\theta_1) = \sigma({{\cdot}},\theta_{1,0})$
if and only if $\theta=\theta_0$. \label{hip4a}
\end{enumerate}

\subsection*{The contrast function}
The main idea of this paper is to fix a new threshold $\tau'$ $(>\tau) 
$ as follows.
We fix a number $\alpha\in(0, 1/2)$ and take
a sequence $\tau_n$ $(>\tau)$ such that
$h_n^\alpha / (\tau_n - \tau) =O(1)$; for example,
$\tau_n=\tau+h_n^\alpha$.
We use $\tau'$ instead of $\tau_n$.
Notice that $\tau' \to \tau$ slowly. 
%
Thus, we introduce the following contrast functions
\begin{eqnarray}
g_n(\theta_1) &=& \sum_{i=1}^n g(i, i-1; \theta_1) \chi_{\{X_{i-1}> 
\tau', X_i>\tau\}}, \label{eq:truecontr0} \\
\ell_n({\theta}) &=& \sum_{i=1}^n \ell(i, i-1; {{{\theta}}}) \chi_{\{X_{i-1}>\tau', X_i>\tau\}}, \label{eq:truecontr}
\end{eqnarray}
where $\chi$ is the indicator function and
\begin{eqnarray*}
g(i, i-1; \theta_1) &=& \log \sigma_{i-1}^2 + \frac{(\Delta_i X)^2} 
{\sigma_{i-1}^2 h_n}, \\
\ell(i, i-1; \theta) &=& \log \sigma_{i-1}^2  + \frac{\left(\Delta_i  
X - b_{i-1}h_n\right)^2}{\sigma_{i-1}^2 h_n}.
\end{eqnarray*}


%


\section{Consistent and asymptotically normal estimators}\label 
{sec:main}
As in Yoshida (1992), we first estimate the parameter belonging to  
the diffusion coefficient,
i.e. $\theta_1$, because, as usual, the estimator of $\theta_1$
has a faster rate of convergence than one of $\theta_2$.
Let $\hat\theta_{1,n}$ denote an estimator of $\theta_1$ satisfying
\begin{equation}
g_n(\hat\theta_{1,n}) = { \inf_{\theta_1} g_n(\theta_1)}.
\label{hat1}
\end{equation}
{{The measurable selection theorem ensures the existence
of such a measurable mapping. }}

\begin{theorem} \label{thm:Theta1}
Under assumptions A\ref{hip1}-A\ref{hip4a}, 
$$
\hat\theta_{1,n}\overset{p}{\to} \theta_{1,0}.
$$
\end{theorem}

\noindent

We consider an estimator $\hat\theta_{2,n}$
of $\theta_2$ that satisfies
\begin{equation}
\ell_n(\hat\theta_{1,n},
\hat\theta_{2,n}) ={\inf_{\theta_2} \ell_n}(\hat\theta_{1,n},\theta_2).
\label{hat2}
\end{equation}

\begin{theorem}\label{thm:Theta2}
Under assumptions A\ref{hip1}-A\ref{hip4a}, 
$$
\hat\theta_{2,n} \overset{p}{\to} \theta_{2,0}.
$$
\end{theorem}

Let
$$
\Sigma\!=\!
\begin{pmatrix}
{ 2 \int \left( \frac{\dA\sigma(x,\theta_{1,0})}{\sigma(x,\theta_ 
{1,0})} \right)^2 \chi_{\{x>\tau\}}\nu_{\theta_0}(\de x)}
& 0\\
0 &
\!\!\!\!\!\!\!\!\!\!\!\!\ { \int \left(\frac{\dB b(x,\theta_{2,0})} 
{\sigma(x,\theta_{1,0})}\right)^2 \chi_{\{x>\tau\}} \nu_{\theta_0} 
(\de x)}
\end{pmatrix}.
$$
Next theorem is the main result in this paper.
\begin{theorem}\label{th0}
{Suppose that the assumptions}
A\ref{hip1}-A\ref{hip4a} 
{{are satisfied.}}
{If $\Sigma$ is non-singular,} then
$$\begin{pmatrix}
\sqrt{n}(\hat\theta_{1,n} - \theta_{1,0})\\
\sqrt{n h_n}(\hat\theta_{2,n} - \theta_{2,0})\\
\end{pmatrix} \overset{d}{\to} N(0, \Sigma^{{-1}}).
$$
\end{theorem}

\section{Proofs} \label{sec:proofs}


\begin{proof}[Proof of Theorem \ref{thm:Theta1}]

First, we will show that
\begin{equation} 
\sup_{\theta_1}
\left| \frac{1}{n} g_n(\theta_1) -G(\theta_1) \right| \overset{p} 
{\to} 0,
\label{uchida:cont1}
\end{equation}
where
$$
G(\theta_1) = \int_{\bf R} \left\{\log \sigma^2(x,\theta_1) + \frac 
{\sigma^2(x,\theta_{1,0})}{\sigma^2(x,\theta_1)}\right\} \chi_{\{x> 
\tau\}} \nu_{\theta_0}(dx).
$$
Noting that
\begin{equation} \label{yoshida:ess1}
\chi_{\{X_{i-1}>\tau', X_i>\tau\}} - \chi_{\{X_{i-1}>\tau'\}} = -  
\chi_{\{X_{i-1}>\tau', X_i\leq\tau\}},
\end{equation}
one has
\begin{eqnarray}
\frac{1}{n}g_n(\theta_1) &=& \frac{1}{n} \sum_{i=1}^n g(i, i-1;  
\theta_1) \chi_{\{X_{i-1}>\tau'\}} \label{uchida:cont1-a} \\
& & - \frac{1}{n}\sum_{i=1}^n g(i, i-1; \theta_1)  \chi_{\{X_{i-1}> 
\tau', X_i\leq\tau\}}. \label{uchida:cont1-b}
\end{eqnarray}
In order to show the uniform convergence of (\ref{uchida:cont1-b}) to  
zero, we estimate
\begin{eqnarray*}
& & \E \left\{ \sup_{\theta_1} \left|\frac{1}{n} \sum_{i=1}^n  g(i,  
i-1; \theta_1)  \chi_{\{X_{i-1}>\tau', X_i\leq\tau\}} \right| \right 
\} \\
&\leq&
\frac{1}{n} \sum_{i=1}^n \left|\left| \ \sup_{\theta_1} \left| \log  
\sigma_{i-1}^2 + \frac{(\Delta_i X)^2}{h_n \sigma_{i-1}^2}\right|\,\,
\right|\right|_p
  P(X_{i-1}>\tau', X_i\leq\tau)^\frac{1}{q}
\end{eqnarray*}
with $1/p + 1/q =1$.
Since A\ref{hip4} and A\ref{hip6} imply that
$$
\sup_{\theta_1} { \left|  \log \sigma_{i-1}^2\right|} \leq \max(|\log  
(K_4)^2|, \sup_{\theta_1}|\sigma_{i-1}^2|)\leq K_4' + C^2(1+|X_{i-1}|) 
^{2C},
$$
it follows from A\ref{hip3} that
$$
  \left|\left| \sup_{\theta_1} \left| \log \sigma_{i-1}^2\right|\,\,  
\right|\right|_p < \infty.
$$
By A\ref{hip4} and the estimate that $\E |X_i - X_{i-1}|^{2p} \leq C  
h_n^{p}$,  
\begin{equation*}
\begin{aligned}
\left|\left|  \sup_{\theta_1} \frac{(\Delta_i X )^2}{h_n \sigma_{i-1} 
^2} \right|\right|_p^{{p}}
&\leq
K \left|\left|  \frac{(\Delta_i X)^2}{h_n}  \right|\right|_p^{{p}} =O(1).
\end{aligned}
\end{equation*}
Moreover, for $k >0$,
\begin{eqnarray}
\sup_i P(X_{i-1}>\tau', X_i\leq\tau) &\leq& \sup_i  P( |X_{i-1}-X_i|  
\geq \tau'-\tau) \nonumber \\
&\leq&  \left(\frac{1}{\tau' - \tau}\right)^k \sup_i  \E |X_{i-1} -  
X_i|^k \nonumber \\
&=& C \left( \frac{h_n^\alpha}{\tau'-\tau} \right)^k (h_n^{1/2- 
\alpha})^k
\label{yoshida:ess2} \\
&=& O \left( h_n^{(1/2-\alpha)k} \right) \to 0 \nonumber
\end{eqnarray}
because 
$h_n^\alpha/(\tau' - \tau) = O(1)$ for $\alpha \in (0,1/2)$.
Thus, we obtain
$$
\sup_{\theta_1} \left|\frac{1}{n} \sum_{i=1}^n  g(i, i-1; \theta)   
\chi_{\{X_{i-1}>\tau', X_i\leq\tau\}} \right| =o_p(1).
$$
In order to prove the uniform convergence of (\ref{uchida:cont1-a})  
to $G$,
it is enough to show that
\begin{equation} \label{cont1-b:prob}
\frac{1}{n} \sum_{i=1}^n g(i, i-1; \theta_1) \chi_{\{X_{i-1}>\tau'\}}  
\overset{p}{\to} G(\theta_1)
\end{equation}
for each $\theta_1$, and
\begin{equation} \label{cont1-b:tight}
\sup_n E\left[ \sup_{\theta_1} \left| \frac{1}{n} \sum_{i=1}^n \delta_ 
{\theta_1} g(i, i-1; \theta_1) \right| \right] < \infty.
\end{equation}
For details, see the proof of Theorem 4.1 in Yoshida (1990).
As in the proof of the uniform convergence of (\ref{uchida:cont1-b}),
we can obtain (\ref{cont1-b:tight}). 
For the proof of (\ref{cont1-b:prob}), we will prove
\begin{eqnarray}
\frac{1}{n} \sum_{i=1}^n \log\sigma^2_{i-1} \chi_{\{X_{i-1} >\tau'\}}
&\overset{p}{\to}& \int_{\bf R} \log\sigma^2(x,\theta) \chi_{\{x>\tau 
\}} \nu_{\theta_0}(dx), \label{uchida:ergod} \\
\frac{1}{n h_n} \sum_{i=1}^n \frac{(\Delta_i X)^2}{\sigma^2_{i-1}}  
\chi_{\{X_{i-1} >\tau'\}}
&\overset{p}{\to}& \int_{\bf R} \frac{\sigma^2(x,\theta_{1,0})} 
{\sigma^2(x,\theta_1)} \chi_{\{x>\tau\}} \nu_{\theta_0}(dx). \label 
{uchida:ergod2}
\end{eqnarray}
For the proof of (\ref{uchida:ergod}),
we set $I_i =  \int_{t_{i-1}}^{t_i} \log\sigma^2_{i-1} \chi_{\{X_ 
{i-1} >\tau'\}} ds$ for $i=1,\ldots,n$.
Note that
$$
\chi_{ \{ X_{i-1} >\tau'\} } \geq \chi_{ \{ \inf_{s \in (t_ 
{i-1},t_i]} X_s >\tau \} }
-\chi_{ \{ \tau < \inf_{s \in (t_{i-1},t_i]} X_s \leq \tau' \} }.
$$
We first estimate $I_i$ for the case that $\log\sigma^2_{i-1} \geq 0 
$. Let $J_i = \chi_{ \{ \log\sigma^2_{i-1} \geq 0 \} }$ for $i=1, 
\ldots,n$.
\begin{eqnarray*}
I_i J_i
&\geq& J_i \int_{t_{i-1}}^{t_i} \log\sigma^2_{i-1} \chi_{\{X_{i-1} > 
\tau'\}}^2 ds \nonumber \\
&\geq& J_i \int_{t_{i-1}}^{t_i} \log\sigma^2_{i-1} \chi_{\{X_{i-1} > 
\tau'\}} \left[ \chi_{ \{ \inf_{s \in (t_{i-1},t_i]} X_s >\tau \} }
-\chi_{ \{ \tau < \inf_{s \in (t_{i-1},t_i]} X_s \leq \tau' \} } 
\right]  ds \nonumber \\
&\geq&
- J_i \int_{t_{i-1}}^{t_i} \log\sigma^2_{i-1} \chi_{\{X_{i-1} > 
\tau'\}} \chi_{ \{ \tau < \inf_{s \in (t_{i-1},t_i]} X_s \leq \tau'  
\} }  ds.
\nonumber \\
& & + J_i \int_{t_{i-1}}^{t_i} \log\sigma^2_{i-1} \chi_{\{X_{i-1} > 
\tau'\}}
\chi_{ \{ \inf_{s \in (t_{i-1},t_i]} X_s >\tau \} }
\chi_{ \{ X_s >\tau \} } ds \nonumber \\ 
&=&
- J_i \int_{t_{i-1}}^{t_i} \log\sigma^2_{i-1} \chi_{\{X_{i-1} > 
\tau'\}} \chi_{ \{ \tau < \inf_{s \in (t_{i-1},t_i]} X_s \leq \tau'  
\} }  ds.
\nonumber \\
& & + J_i \int_{t_{i-1}}^{t_i} \log\sigma^2_{i-1} \chi_{\{X_{i-1} >\tau'\}}
\left[ \chi_{ \{ \inf_{s \in (t_{i-1},t_i]} X_s >\tau \} }-1 \right]
\chi_{ \{ X_s >\tau \} } ds
\nonumber \\
& & + J_i \int_{t_{i-1}}^{t_i} \log\sigma^2_{i-1}
\left[ \chi_{\{X_{i-1} >\tau'\}}  -1 \right]
\chi_{ \{ X_s >\tau \} } ds
\nonumber \\
& & + J_i \int_{t_{i-1}}^{t_i} \log\sigma^2_{i-1}  \chi_{ \{ X_s >  
\tau \} } ds. \nonumber
\end{eqnarray*}
Hence,
\begin{eqnarray}
J_i \left( I_i - \int_{t_{i-1}}^{t_i} \log\sigma^2(X_s, \theta_1)  
\chi_{\{X_s >\tau \}}  ds \right) &\geq& \Xi_i^{(1)}, \label{lowbound}
\end{eqnarray}
where
\begin{eqnarray}
\Xi_i^{(1)}
&=&
- J_i \int_{t_{i-1}}^{t_i} \log\sigma^2_{i-1} \chi_{\{X_{i-1} > 
\tau'\}} \chi_{ \{ \tau < \inf_{s \in (t_{i-1},t_i]} X_s \leq \tau'  
\} }  ds.
\label{erg2} \\
& & - J_i \int_{t_{i-1}}^{t_i} \log\sigma^2_{i-1} \chi_{\{X_{i-1} >  
\tau'\}} \chi_{ \{ \inf_{s \in (t_{i-1},t_i]} X_s \leq \tau \} }
\chi_{ \{ X_s >\tau \} }ds
\label{erg3} \\
& & - J_i \int_{t_{i-1}}^{t_i} \log\sigma^2_{i-1} \chi_{\{X_{i-1}  
\leq \tau'\}}
\chi_{ \{ X_s >\tau \} }ds  \label{erg4} \\
& & + J_i \int_{t_{i-1}}^{t_i} \left\{ \log\sigma^2_{i-1} - \log 
\sigma^2(X_s, \theta_1) \right\} \chi_{ \{ X_s >\tau \} } ds. \label 
{erg5}
\end{eqnarray}
Next, noting that
\begin{eqnarray*}
I_i J_i &=& J_i \int_{t_{i-1}}^{t_i} \log\sigma^2_{i-1} \chi_{\{X_ 
{i-1} >\tau'\}}
\left[ \chi_{ \{ \inf_{s \in (t_{i-1},t_i]} X_s >\tau \} } +\chi_{ \{ \inf_{s \in (t_{i-1},t_i]} X_s  \leq \tau \} } \right]ds \\
&=& J_i \int_{t_{i-1}}^{t_i} \log\sigma^2_{i-1} \chi_{\{X_{i-1} > 
\tau'\}}
\chi_{ \{ \inf_{s \in (t_{i-1},t_i]} X_s >\tau \} }ds \\
& & + J_i \int_{t_{i-1}}^{t_i} \log\sigma^2_{i-1} \chi_{\{X_{i-1} > 
\tau'\}}
\chi_{ \{ \inf_{s \in (t_{i-1},t_i]} X_s  \leq \tau \} } ds \\
&\leq& J_i \int_{t_{i-1}}^{t_i} \log\sigma^2_{i-1} \chi_{ \{ X_s > 
\tau \} }ds
  + J_i \int_{t_{i-1}}^{t_i} \log\sigma^2_{i-1} \chi_{\{X_{i-1} > 
\tau'\}}
\chi_{ \{ \inf_{s \in (t_{i-1},t_i]} X_s  \leq \tau \} } ds,
\end{eqnarray*}
we obtain that
\begin{eqnarray}
J_i \left( I_i - \int_{t_{i-1}}^{t_i} \log\sigma^2(X_s, \theta_1)  
\chi_{\{X_s >\tau \}}  ds \right) &\leq& \Xi_i^{(2)}, \label{upbound}
\end{eqnarray}
where
\begin{eqnarray*}
\Xi_i^{(2)}
&=& \int_{t_{i-1}}^{t_i} \left\{ \log\sigma^2_{i-1} - \log\sigma^2 
(X_s, \theta_1) \right\}  \chi_{ \{ X_s >\tau \} }ds \\
& &
  + J_i \int_{t_{i-1}}^{t_i} \log\sigma^2_{i-1} \chi_{\{X_{i-1} > 
\tau'\}}
\chi_{ \{ \inf_{s \in (t_{i-1},t_i]} X_s  \leq \tau \} } ds.
\end{eqnarray*}
It follows from (\ref{lowbound}) and (\ref{upbound}) that
$$
\left| J_i \left( I_i - \int_{t_{i-1}}^{t_i} \log\sigma^2(X_s,  
\theta_1) \chi_{\{X_s >\tau \}}  ds \right) \right|
\leq \max \{ |\Xi_i^{(1)}|,  |\Xi_i^{(2)}| \}.
$$
For the estimate of (\ref{erg2}), we set $\tilde{\tau} = \tau' + h_n^ 
\alpha$, where $\alpha \in (0,1/2)$.
\begin{eqnarray*}
& & 
E \left[ \left|
J_i \int_{t_{i-1}}^{t_i} \log\sigma^2_{i-1} \chi_{\{X_{i-1} >  
\tau'\}} \chi_{ \{ \tau <  \inf_{s \in (t_{i-1},t_i]} X_s \leq \tau'  
\} } ds
\right| \right] \\
&\leq& 
E \left[ \left|
\int_{t_{i-1}}^{t_i} \log\sigma^2_{i-1} \left\{ \chi_{\{X_{i-1} >  
\tau'\}} -  \chi_{ \{ X_{i-1} > \tilde{\tau} \}} \right\}
\chi_{ \{ \tau < \inf_{s \in (t_{i-1},t_i]} X_s \leq \tau' \} } ds
\right| \right] \\
& & + 
E \left[ \left|
\int_{t_{i-1}}^{t_i} \log\sigma^2_{i-1} \chi_{ \{ X_{i-1} > \tilde 
{\tau} \}}
\chi_{ \{ \tau < \inf_{s \in (t_{i-1},t_i]} X_s \leq \tau' \} } ds
\right| \right] \\
&\leq& 
E \left[ \left|
\int_{t_{i-1}}^{t_i} \log\sigma^2_{i-1}  \chi_{ \{ \tau' < X_{i-1}  
\leq \tilde{\tau} \}}  ds
\right| \right] \\
& & + 
E \left[ \left|
\int_{t_{i-1}}^{t_i} \log\sigma^2_{i-1} \chi_{ \{ X_{i-1} > \tilde 
{\tau} \}}
\chi_{ \{ \inf_{s \in (t_{i-1},t_i]} X_s \leq \tau' \} } ds
\right| \right] \\
&\leq& h_n C \left[   P(\tau' < X_{i-1} \leq \tilde{\tau})^{1/2}
+  P( \sup_{s \in (t_{i-1},t_i]}  | X_{i-1} - X_s | > h_n^\alpha )^ 
{1/2}  \right]
=o(h_n).
\end{eqnarray*}
Concerning the estimate of (\ref{erg3}), 
\begin{eqnarray*}
& &  E \left[ \left|
J_i \int_{t_{i-1}}^{t_i} \log\sigma^2_{i-1} \chi_{\{X_{i-1} > \tau'\}}
\chi_{ \{ \inf_{s \in (t_{i-1},t_i]} X_s \leq \tau \} } \chi_{\{X_s > 
\tau \}}  ds
\right| \right] \\
&\leq&  h_n C   P( \sup_{s \in (t_{i-1},t_i]} | X_{i-1}  -X_s | > h_n^ 
\alpha )^{1/2}
=o(h_n).
\end{eqnarray*}
In order to estimate (\ref{erg4}), we set $\tau'' = \tau - h_n^\alpha 
$, where $\alpha \in (0,1/2)$.
\begin{eqnarray*}
& &  E \left[ \left|
J_i \int_{t_{i-1}}^{t_i} \log\sigma^2_{i-1} \chi_{\{X_{i-1} \leq  
\tau'\}}
  \chi_{\{X_s >\tau \}}  ds
\right| \right] \\
&\leq&  E \left[ \left|
\int_{t_{i-1}}^{t_i} \log\sigma^2_{i-1} \left\{ \chi_{\{X_{i-1} \leq  
\tau'\}} -  \chi_{ \{\tau''< X_{i-1} \leq \tau'\}} \right\}
\chi_{\{X_s >\tau \}} ds
\right| \right] \\
& & +
E \left[ \left|
\int_{t_{i-1}}^{t_i} \log\sigma^2_{i-1} \chi_{ \{\tau''< X_{i-1} \leq  
\tau'\}} \chi_{\{X_s >\tau \}} ds
\right| \right] \\
&\leq&
E \left[ \left|
\int_{t_{i-1}}^{t_i} \log\sigma^2_{i-1}  \chi_{\{X_{i-1} \leq \tau''\}}
\chi_{\{X_s >\tau \}} ds
\right| \right]
+
E \left[ \left|
\int_{t_{i-1}}^{t_i} \log\sigma^2_{i-1} \chi_{ \{\tau''< X_{i-1} \leq  
\tau'\}} ds
\right| \right] \\
&\leq& h_n C \left[  P( \sup_{s \in (t_{i-1},t_i]} | X_s  -X_{i-1} |  
 > h_n^\alpha )^{1/2} + P(\tau'' < X_{i-1} \leq \tau')^{1/2} \right]
=o(h_n).
\end{eqnarray*}
As for the estimate of (\ref{erg5}),
\begin{eqnarray*}
E \left[ \left|
J_i \int_{t_{i-1}}^{t_i} \left\{ \log\sigma^2_{i-1} - \log\sigma^2 
(X_s, \theta_1) \right\}   ds
\right| \right]
&\leq& C h_n^{3/2} = o(h_n).
\end{eqnarray*}
Thus, we obtain
\begin{equation} \label{ergod1}
E \left[
| \Xi_i^{(1)} |
\right] = o(h_n).
\end{equation}
Moreover,
\begin{eqnarray}
E \left[
| \Xi_i^{(2)} |
\right]
&\leq& 
E \left[ \left|
\int_{t_{i-1}}^{t_i} \left\{ \log\sigma^2_{i-1} - \log\sigma^2(X_s,  
\theta_1) \right\}  \chi_{ \{ X_s >\tau \} }ds
\right| \right]  \nonumber \\
& & + 
E \left[ \left| \int_{t_{i-1}}^{t_i} \log\sigma^2_{i-1} \chi_{\{X_ 
{i-1} >\tau'\}}
\chi_{ \{ \inf_{s \in (t_{i-1},t_i]} X_s  \leq \tau \} } ds \right|  
\right] \nonumber \\
&=& 
o(h_n). \label{ergod2}
\end{eqnarray}
It follows from (\ref{ergod1}) and (\ref{ergod2}) that
\begin{equation} \label{erg10}
E \left[ \left|
J_i \int_{t_{i-1}}^{t_i}
\left\{ \log\sigma^2_{i-1} \chi_{\{X_{i-1} >\tau'\}}
- \log \sigma^2(X_s,\theta_1) \chi_{\{ X_s > \tau \}} \right\}  ds
\right| \right] =o(h_n).
\end{equation}
For the case that $\log\sigma^2_{i-1} < 0$,
in a similar way as above,
we can show that
\begin{equation} \label{erg10-b}
E \left[ \left|
(1-J_i) \int_{t_{i-1}}^{t_i}
\left\{ \log\sigma^2_{i-1} \chi_{\{X_{i-1} >\tau'\}}
- \log \sigma^2(X_s,\theta_1) \chi_{\{ X_s > \tau \}} \right\}  ds
\right| \right] =o(h_n).
\end{equation}
Therefore, we have
$$
E \left[ \left|
\int_{t_{i-1}}^{t_i}
\left\{ \log\sigma^2_{i-1} \chi_{\{X_{i-1} >\tau'\}}
- \log \sigma^2(X_s,\theta_1) \chi_{\{ X_s > \tau \}} \right\}  ds
\right| \right] =o(h_n)
$$
and consequently,
\begin{equation} \label{erg11}
\left|
\frac{1}{n h_n} \sum_{i=1}^n
\int_{t_{i-1}}^{t_i}
\left\{ \log\sigma^2_{i-1} \chi_{\{X_{i-1} >\tau'\}}
- \log \sigma^2(X_s,\theta_1) \chi_{\{ X_s > \tau \}} \right\}  ds
\right| =o_p(1).
\end{equation}
Moreover, by the ergodic theorem,
$$
\frac{1}{n h_n } \int_0^{nh_n} \log \sigma^2(X_s,\theta_1) \chi_{\{ X_s > \tau \}}  ds
\overset{p}{\to}
\int_{\bf R} \log \sigma^2(x,\theta_1) \chi_{\{ x > \tau \}} \nu_ 
{\theta_0}(dx),
$$
which completes the proof of (\ref{uchida:ergod}).
For the proof of (\ref{uchida:ergod2}),
we set
$$\Xi_i = \frac{1}{n h_n} \frac{(\Delta_i X)^2}{\sigma^2_{i-1}} \chi_ 
{\{X_{i-1} >\tau'\}}.
$$
By Lemma 9 of Genon-Catalot and Jacod (1993), it is enough to show that
\begin{eqnarray}
\sum_{i=1}^n \E_{\theta_0} \left\{ \Xi_i | \mathcal F_{i-1}\right\} & 
\overset{p}{\to}&
\int_{\bf R} \frac{\sigma^2(x,\theta_{1,0})}{\sigma^2(x,\theta_1)}  
\chi_{\{x>\tau\}} \nu_{\theta_0}(dx), \label{eq:Xi1}\\
\sum_{i=1}^n \E_{\theta_0}\left\{ \left(\Xi_i \right)^2 | \mathcal F_ 
{i-1}\right\} &\overset{p}{\to}& 0,  \label{eq:Xi2}
\end{eqnarray}
where $\mathcal F_{i-1}$ denotes the {\it history} up to the time $t_ 
{i-1}$.
In order to evaluate
$
\E_{\theta_0}\left\{ \left(\Delta_i X \right)^2|\mathcal F_{i-1}\right 
\},
$
we can use a well known It\^o-Taylor expansion:
\begin{eqnarray*}
\E_{\theta_0}(\phi(X_i,X_{i-1})  |\mathcal F_{i-1}) &=&
  \phi(X_{i-1},X_{i-1})
+ h_n L_{\theta_0} \phi(X_{i-1},X_{i-1})
\\ &&
+
\frac12 h_n^2 L^2_{\theta_0} \phi(X_{i-1},X_{i-1})
\\ &&
+\int_{t_{i-1}}^{t_i}\int_{t_{i-1}}^t
E_{\theta_0}\Bigl\{L^2_{\theta_0}\phi(X(s),X_{i-1})
\\ &&
-L^2_{\theta_0}\phi(X_{i-1},X_{i-1})|\mathcal F_{i-1}\Bigr\}\>ds\>dt
\end{eqnarray*}
for appropriate functions $\phi(x,y)$,
where
$L_\theta \phi(x,y)= \frac12 \sigma^2(x,\theta_1)
\frac{\partial^2}{\partial x^2} \phi(x,y) + b(x,\theta_2) \frac 
{\partial}{\partial x}\phi(x,y)$.
Hence
\begin{equation} \label{I-T1}
\E_{\theta_0}\left\{ \left(\Delta_i X \right)^2|\mathcal F_{i-1}\right\}
= h_n \sigma^{*2}_{i-1}+ R(h_n^2,  X_{i-1}),
\end{equation}
where $R(\cdot, \cdot)$ is defined in \eqref{eq:Rn}.
Thus
\begin{eqnarray*}
\sum_{i=1}^n E_{\theta_0} \left\{ \Xi_{i} | \mathcal F_{i-1} \right\}
&=&  \frac{1}{n} \sum_{i=1}^n  \frac{ (\sigma^*_{i-1})^2 }{\sigma^2_ 
{i-1}} \chi_{ \{X_{i-1} >\tau'\} }
+ \frac{h_n}{n} \sum_{i=1}^n  R(1,X_{i-1}) \\
&\overset{p}{\to}&
\int_{\bf R} \frac{\sigma^2(x,\theta_{1,0})}{\sigma^2(x,\theta_1)}  
\chi_{\{x>\tau\}} \nu_{\theta_0}(dx)
\end{eqnarray*}
and in a similar way, we can show (\ref{eq:Xi2}).
This completes the proof of (\ref{uchida:cont1}).

%

Next, we see that 
$G$ attains to its minimum only at $\theta_{1,0}$
by noting that
$$
\frac{d}{dx} \left(\log x+\frac{a}{x}\right) = \frac{1}{x}-\frac{a} 
{x^2}=\frac{x-a}{x^2}.
$$
Hence, for any $\epsilon >0$, $\inf_{\theta_1: |\theta_1-\theta_ 
{1,0}| \geq \epsilon} G(\theta_1) > G(\theta_{1,0})$.
This implies that if $| \theta_1 -\theta_{1,0} | \geq \epsilon$, then  
$G(\theta_1) > G(\theta_{1,0}) + \eta$ for some $\eta >0$.
Therefore,
\begin{eqnarray}
P \left( | \hat{\theta}_{1,n} - \theta_{1,0} | \geq \epsilon \right)
&\leq& P \left( G(\hat{\theta}_{1,n}) > G(\theta_{1,0}) + \eta  
\right) \nonumber \\
&\leq& 2 P \left( \sup_{\theta_1} \left| \frac{1}{n}g_n(\theta_1) -G 
(\theta_1) \right| > {\eta}/{3} \right) \label{uchida:sup} \\
& & + P \left( \frac{1}{n}g_n(\hat{\theta}_{1,n})- \frac{1}{n}g_n 
(\theta_{1,0}) > {\eta}/{3} \right). \nonumber
\end{eqnarray}
By using (\ref{uchida:cont1}), the probability of (\ref{uchida:sup})  
converges 0.
Furthermore, it follows from
(\ref{hat1})
that
$$
P \left( \frac{1}{n}g_n(\hat{\theta}_{1,n})- \frac{1}{n}g_n(\theta_ 
{1,0}) > {\eta}/{3} \right)
\leq P \left( \frac{1}{n}g_n(\hat{\theta}_{1,n}) > \frac{1}{n}g_n 
(\theta_{1,0})  \right)
\rightarrow 0.
$$
This competes the proof.
\end{proof}

\vspace{10pt}

\begin{proof}[Proof of Theorem \ref{thm:Theta2}]
We need to prove that
\begin{equation} \label{convT1}
\sup_{\theta_2}\left|
\frac{1}{n h_n}\left(
\ell_n(\hat\theta_{1,n}, \theta_2) - \ell_n(\hat\theta_{1,n}, \theta_ 
{2,0})
\right) - L(\theta_2)\right| \overset{p}{\to} 0,
\end{equation}
where
\begin{equation*}
L(\theta_2) =  \int_{\bf R} \left(\frac{b(x,\theta_2)-b(x,\theta_ 
{2,0})}{\sigma(x,\theta_{1,0})}\right)^2\chi_{\{x>\tau\}}\nu_ 
{\theta_0}(dx).
\label{limT2}
\end{equation*}
An easy computation together with (\ref{yoshida:ess1}) yields that
$$
\frac{1}{n h_n}\left(
\ell_n(\hat\theta_{1,n}, \theta_2) - \ell_n(\hat\theta_{1,n}, \theta_ 
{2,0})
\right)
=   \psi_{1,n}(\theta_2) + \psi_{2,n}(\theta_2) + \psi_{3,n} 
(\theta_2) + R_n(\theta_2),
$$
where $\hat\sigma_i = \sigma(X_i, \hat \theta_{1,n})$,
\begin{eqnarray*}
\psi_{1,n}(\theta_2) &=& -\frac{2}{n h_n}\sum_{i=1}^n \frac{(b_{i-1}- 
b_{i-1}^*)
\int_{t_{i-1}}^{t_i}\sigma(X_s, \theta_{1,0})dW_s}{\hat\sigma_{i-1} 
^2 }\chi_{\{X_{i-1}>\tau'\}}, \label{eq:A} \\
\psi_{2,n}(\theta_2) &=& -\frac{2}{n h_n}\sum_{i=1}^n \frac{(b_{i-1}- 
b_{i-1}^*)
\int_{t_{i-1}}^{t_i}b(X_s, \theta_{2,0})ds}{\hat\sigma_{i-1}^2 }\chi_ 
{\{X_{i-1}>\tau'\}}, \label{eq:B} \\
\psi_{3,n}(\theta_2) &=& \frac{1}{n} \sum_{i=1}^n\frac{b_{i-1}^2-b_ 
{i-1}^{*2}}{\hat\sigma^2_{i-1}} \chi_{\{X_{i-1}>\tau'\}},
\label{eq:C} \\
R_n(\theta_2) &=&
\frac{1}{n h_n}\sum_{i=1}^n
\biggl\{
  \frac{(\Delta_i X - b_{i-1} h_n)^2}{\hat\sigma_{i-1}^2 h_n}
-\frac{(\Delta_i X-b_{i-1}^* h_n)^2}{\hat\sigma_{i-1}^2 h_n}\biggr\} 
\chi_{\{X_{i-1}>\tau', X_i \leq \tau \}}.
\end{eqnarray*}
We first estimate $R_n(\theta_2)$.
\begin{eqnarray*}
E \left[ \sup_{\theta_2} |R_n(\theta_2) | \right]
&\leq& \frac{1}{n h_n } \sum_{i=1}^n
E \left[ \sup_{\theta_2} \left| \frac{(\Delta_i X - b_{i-1} h_n)^2 -  
(\Delta_i X-b_{i-1}^* h_n)^2}{\hat\sigma_{i-1}^2 h_n}
\right|^2  \right]^{1/2} \\
& & \times
P( X_{i-1}>\tau', X_i \leq \tau )^{1/2}  \\
&\leq& \frac{1}{h_n^{1/2}} \times C \left( \frac{h_n^\alpha}{\tau'- 
\tau} \right)^{k/2} (h_n^{1/2-\alpha})^{k/2} \\
&=& O \left( h_n^{k/4- \alpha k /2 -1/2} \right) \rightarrow 0,
\end{eqnarray*}
where we took $k > 2/(1- 2 \alpha)$ in (\ref{yoshida:ess2}).
This yields that $\sup_{\theta_2} | R_n(\theta_2) | =o_p(1)$.
Next, $\psi_{2,n}(\theta_2)$ can be rewritten as 
\begin{eqnarray*}
\psi_{2,n}(\theta_2) &=& \psi_{2,n}^{(1)}(\theta_2) + \psi_{2,n}^{(2)} 
(\theta_2) + \psi_{2,n}^{(3)}(\theta_2),
\end{eqnarray*}
where
\begin{eqnarray*}
\psi_{2,n}^{(1)}(\theta_2) &=& -\frac{2}{n }\sum_{i=1}^n \frac{(b_ 
{i-1}-b_{i-1}^*)
b_{i-1}^* }{\sigma_{i-1}^{*2} }\chi_{\{X_{i-1}>\tau'\}}, \\
\psi_{2,n}^{(2)}(\theta_2) &=& -\frac{2}{n h_n}\sum_{i=1}^n \frac{(b_ 
{i-1}-b_{i-1}^*)
\int_{t_{i-1}}^{t_i} \left\{ b(X_s, \theta_{2,0}) -b_{i-1}^* \right\}  
ds}{\sigma_{i-1}^{*2} }\chi_{\{X_{i-1}>\tau'\}}, \\
\psi_{2,n}^{(3)}(\theta_2) &=& -\frac{2}{n h_n}\sum_{i=1}^n (b_{i-1}- 
b_{i-1}^*)\int_{t_{i-1}}^{t_i}b(X_s, \theta_{2,0})ds
\left(\frac{1}{\hat\sigma_{i-1}^2} - \frac{1}{\sigma_{i-1}^{*2}} 
\right)\chi_{\{X_{i-1}>\tau'\}}.
(\theta_2).
\end{eqnarray*}
By noting that for $p, K >0$,
\begin{eqnarray*}
\left| \left| \int_{t_{i-1}}^{t_i} \left\{ b(X_s,\theta_{2,0}) -b_ 
{i-1}^* \right\} ds \right| \right|_p
&\leq& C h_n^{3/2}, \\
\left|\frac{1}{\hat\sigma_{i-1}^2} - \frac{1}{\sigma_{i-1}^{*2}} 
\right| \leq C\left|
\sigma_{i-1}^{*2} - \hat\sigma_{i-1}^2 \right| &\leq& |\hat\theta_ 
{1,n} - \theta_{1,0}| K(1+|X_{i-1}|)^K, 
\end{eqnarray*}
one has that for $p,q>0$ with $1/p +1/q =1$,
\begin{eqnarray*}
E \left[ \sup_{\theta_2} \left| \psi_{2,n}^{(2)}(\theta_2) \right|  
\right]
&\leq& \frac{1}{n h_n} \sum_{i=1}^n \left| \left| \ \sup_{\theta_2}  
\left|
\frac{b_{i-1}-b_{i-1}^*}{\sigma_{i-1}^{*2}} \right| \ \right| \right|_p
\\
& & \times
\left| \left| \int_{t_{i-1}}^{t_i} \left\{ b(X_s,\theta_{2,0}) -b_ 
{i-1}^* \right\} ds \right| \right|_q \\
&\leq& \frac{C}{n h_n} n h_n^{3/2} \rightarrow 0,
\end{eqnarray*}
and
\begin{eqnarray*}
\sup_{\theta_2} \left| \psi_{2,n}^{(3)}(\theta_2) \right|
&\leq& \left| \hat{\theta}_{1,n} -\theta_{1,0} \right| \frac{K}{n  
h_n} \sum_{i=1}^n \sup_{\theta_2} \left| b_{i-1} - b_{i-1}^* \right|
\\
& & \times
\left| \int_{t_{i-1}}^{t_i} b(X_s,\theta_{2,0})  ds \right|
(1+|X_{t_{i-1}}|)^K \\
&=& o_p(1) \times O_p(1) = o_p(1).
\end{eqnarray*}
As in the proof of the uniform convergence of (\ref{uchida:cont1-a}),
\begin{eqnarray*}
\sup_{\theta_2} \left| \psi_{2,n}^{(1)}(\theta_2) +
2\int  \frac{(b(x,\theta_2)-b(x, \theta_{2,0}) )b(x, \theta_{2,0})} 
{\sigma^2(x,\theta_{1,0})}\chi_{\{x>\tau\}}
\nu_{\theta_0} (dx) \right| = o_p(1).
\end{eqnarray*}
Furthermore, since one estimates
\begin{eqnarray*}
& & \sup_{\theta_2} \left| \psi_{3,n}(\theta_2) -
\frac{1}{n} \sum_{i=1}^n \frac{b_{i-1}^2-b_{i-1}^{*2}}{\sigma^{*2}_ 
{i-1}} \chi_{\{X_{i-1}>\tau'\}} \right| \\
&\leq& \left| \hat{\theta}_{1,n} -\theta_{1,0} \right| \frac{K}{n  
h_n} \sum_{i=1}^n
\sup_{\theta_2} \left| b_{i-1}^2-b_{i-1}^{*2} \right| (1+|X_{t_ 
{i-1}}|)^K \\
&=& o_p(1) \times O_p(1) = o_p(1),
\end{eqnarray*}
we obtain
\begin{equation*}
\sup_{\theta_2} \left| \psi_{3,n}(\theta_2) - \int \frac{b(x,\theta_2) 
^2-b(x,\theta_{2,0})^2}{\sigma^2(x,\theta_{1,0})} \chi_{\{x>\tau\}}
\nu_{\theta_0}(\de x) \right| = o_p(1).
\label{eq:C1}
\end{equation*}
Therefore, we see that
$$
\sup_{\theta_2} \left| \psi_{2,n}(\theta_2) + \psi_{3,n}(\theta_2) - L 
(\theta_2) \right| =o_p(1).
$$
To estimate $\psi_1(\theta_1)$, we consider the following process
$$
M_n(\theta)=\int_0^{n h_n}
\sum_{i=1}^n \frac{(b_{i-1} -b_{i-1}^*) \sigma(X_s, \theta_{1,0})}{n  
h_n \sigma_{i-1}^2 }{\bf 1}_i(s) dW_s,
$$
where
${\bf 1}_i(s) = \chi_{\{X_{i-1}>\tau'\}} \chi_{\{t_{i-1} \leq s \leq  
t_i\}}$.
We will prove the followings:
there exists a constant $\alpha >2$ such that for any $\theta$ and $ 
\theta'$,
\begin{eqnarray}
M_n(\theta) &\overset{p}{\to}& 0,
\label{sup1} \\
E|M_n(\theta)|^\alpha &\leq& C,
\label{tight2} \\
\E\left| M_n(\theta) - M_n(\theta')\right|^\alpha &\leq& C |\theta -  
\theta'|^\alpha,
\label{tight1}
\end{eqnarray}
where $C$ is a constant independent of $\theta$, $\theta'$ and $n$.
If (\ref{sup1})-(\ref{tight1}) are satisfied,
by Theorem 20 in Appendix of Ibragimov and Has'minskii (1981)
or Lemma 3.1 of Yoshida (1990),
we can show that $\sup_\theta |M_n(\theta)| \overset{p}{\to} 0$.
In fact, (\ref{tight2})-(\ref{tight1}) ensure that
the family of distributions of $\{ M_n(\cdot)\}$ on $C(\Theta)$ with  
sup-norm is tight.
Hence, one will be able to prove that
\begin{equation} \label{uchida:psi1}
\sup_{\theta_2}\left|\psi_{1,n}(\theta_2)\right|
=2 \sup_{\theta_2}\left|M_n(\hat\theta_{1,n}, \theta_2)\right|
  \leq 2 \sup_{\theta} |M_n(\theta)| \overset{p}{\to} 0.
\end{equation}
The proof of (\ref{tight1}) is as follows.
Let us define
$$
\begin{aligned}
f_{i-1}(\theta, \theta') &=  \frac{b_{i-1}(\theta_2)-b_{i-1}^*} 
{\sigma_{i-1}^2(\theta_1) }
-\frac{b_{i-1}(\theta_2')-b_{i-1}^*}{\sigma_{i-1}^2(\theta_1') }\\
&=\frac{b_{i-1}(\theta_2)-b_{i-1}(\theta_2')}{\sigma_{i-1}^2 
(\theta_1') }
+ (b_{i-1}(\theta_2)-b_{i-1}^*)\left(
\frac{1}{\sigma_{i-1}^2(\theta_1)}-\frac{1}{\sigma_{i-1}^2(\theta_1')}
\right).
\end{aligned}
$$
By the Burkholder-Davis-Gundy and Jensen inequalities,
$$
\begin{aligned}
& E |M_n(\theta) - M_n(\theta')|^\alpha \\
&=
\frac{1}{(n h_n)^\alpha}\E  \left|
\int_0^{n h_n}  \sum_{i=1}^n f_{i-1}(\theta,\theta')
\sigma(X_s, \theta_{1,0}) {\bf 1}_i(s) \de W_s
\right|^\alpha
\\
&\leq
\frac{C_\alpha}{(n h_n)^\alpha} \E\left(
\sum_{i=1}^n \int_0^{n h_n} \left(f_{i-1}(\theta,\theta')\sigma(X_s,  
\theta_{1,0})\right)^2
  {\bf 1}_i(s) \de s
\right)^\frac{\alpha}{2}\\
&\leq
\frac{C_\alpha}{(n h_n)^\alpha} n^{\alpha/2-1} \sum_{i=1}^n \E\left(
  \int_{t_{i-1}}^{t_i} \left(f_{i-1}(\theta,\theta')\sigma(X_s,  
\theta_{1,0})\right)^2
ds
\right)^\frac{\alpha}{2}\\
&\leq
\frac{C_\alpha}{(n h_n)^\alpha} (n h_n)^{\alpha/2-1} \sum_{i=1}^n \E 
\left(
\int_{t_{i-1}}^{t_i}
\left|
  f_{i-1}(\theta,\theta')\sigma(X_s, \theta_{1,0})\right|^\alpha
ds \right).
\end{aligned}
$$
Moreover, 
it follows from A\ref{hip5}-A\ref{hip6},
$$
|f_{i-1}(\theta,\theta')|^\alpha \leq K (1+|X_{i-1}|)^K |\theta -  
\theta'|^\alpha,
$$
which completes the proof of (\ref{tight1}).
In the similar way, we can show (\ref{tight2}).
For the proof of (\ref{sup1}),
we set $g_{i} = (b_i(\theta_2) - b_i^*)/\sigma_{i}^2(\theta_1)$ and
$$
\begin{aligned}
\E \left( M_n(\theta)\right)^2
\leq& \frac{1}{n^2h_n^2} \sum_{i=1}^n  \E\left\{\int_{t_{i-1}}^{t_i}
   g_{i-1}^2 \sigma^2(X_s, \theta_{1,0})
ds
\right\} \rightarrow 0, 
\end{aligned}
$$
which completes the proof of (\ref{sup1}). Thus, one can show (\ref{uchida:psi1}), which completes the proof of (\ref{convT1}).
Finally, note that
for any $\epsilon >0$, $\inf_{\theta_2: |\theta_2-\theta_{2,0}| \geq  \epsilon} L(\theta_2) > 0$
because $L$ attains to its minimum only at $\theta_{2,0}$.
As in the proof of Theorem \ref{thm:Theta1},
we can show the consistency of $\hat{\theta}_{2,n}$.
This completes the proof.
\end{proof}



\begin{proof}[Proof of Theorem \ref{th0}]
First, we study the asymptotic normality of the score function.
Let
$$
{\cal L}_n =
\left(
\begin{array}{c}
-\frac{1}{\sqrt{n}} \delta_{\theta_1} g_n(\theta_{1,0}) \\
-\frac{1}{\sqrt{n h_n}} \delta_{\theta_2} \ell_n(\hat\theta_{1,n},  
\theta_{2,0})
\end{array}
\right), \quad
\bar{\cal L}_n =
\left(
\begin{array}{c}
-\frac{1}{\sqrt{n}} \delta_{\theta_1} \bar{g}_n(\theta_{1,0}) \\
-\frac{1}{\sqrt{n h_n}} \delta_{\theta_2} \bar{\ell}_n(\theta_0)
\end{array}
\right),
$$
where
\begin{eqnarray*}
\bar{g}_n(\theta_1) &=& \sum_{i=1}^n g(i,i-1;\theta_1) \chi_{ \{ X_ 
{i-1} > \tau' \} }, \\
\bar{\ell}_n(\theta) &=& \sum_{i=1}^n \ell(i,i-1;\theta) \chi_{ \{ X_ 
{i-1} > \tau' \} }.
\end{eqnarray*}
In order to show that
${\cal L}_n - \bar{\cal L}_n = o_p(1)$,
it is sufficient to show that
\begin{eqnarray}
A_n := \frac{1}{\sqrt{n}} \left( \delta_{\theta_1} g_n(\theta_{1,0}) -
  \delta_{\theta_1} \bar{g}_n(\theta_{1,0}) \right)
&=& o_p(1), \label{score1} \\
B_n := \frac{1}{\sqrt{n h_n}} \left(
\delta_{\theta_2} \ell_n(\hat{\theta}_{1,n}, \theta_{2,0}) -
\delta_{\theta_2} \bar\ell_n(\theta_0) \right)
&=& o_p(1). \label{score2}
\end{eqnarray}
For the proof of (\ref{score1}), one estimates
\begin{eqnarray*}
E|A_n| &\leq& \frac{1}{\sqrt{n}} \sum_{i=1}^n
E \left| \delta_{\theta_1} g(i-1,i;\theta_{1,0}) \chi_{\{X_{i-1}> 
\tau', X_i\leq\tau\}} \right| \\
&\leq& \frac{C}{\sqrt{n}} \sum_{i=1}^n
\left| \left| \frac{\dA\sigma_{i-1}^*}{\sigma_{i-1}^*} \left(1 - \frac 
{\left(\Delta_i X\right)^2}{h_n{{\sigma_{i-1}^{*2}}}}
\right) \right| \right|_2 \times  O \left( h_n^{(1/4-\alpha/2)k}  
\right) \\
&\leq& C \sqrt{n} h_n \times O  \left( h_n^{(1/4-\alpha/2)k-1}  
\right) \rightarrow 0,
\end{eqnarray*}
where we took $k > 4/(1-2\alpha)$ in (\ref{yoshida:ess2}).
For the proof of (\ref{score2}), one has that for $\epsilon >0$,
\begin{eqnarray*}
|B_n| \chi_{\{ | \hat{\theta}_{1,n} - \theta_{1,0} | < \epsilon \} } & 
\leq& \frac{1}{\sqrt{n h_n}} \sum_{i=1}^n
\sup_{\theta_1}
\left| \delta_{\theta_1} \delta_{\theta_2} \ell(i,i-1;\theta_1,  
\theta_{2,0}) \right| \left| \hat{\theta}_{1,n} -\theta_{1,0} \right| \\
& &
+ \frac{1}{\sqrt{n h_n}} \sum_{i=1}^n  \left| \delta_{\theta_2} \ell 
(i,i-1;\theta_0) \chi_{\{X_{i-1}>\tau', X_i\leq\tau\}} \right|.
\end{eqnarray*}
As in the proof of (\ref{score1}),
$
\frac{1}{\sqrt{n h_n}} \sum_{i=1}^n  \left| \delta_{\theta_2} \ell 
(i,i-1;\theta_0) \chi_{\{X_{i-1}>\tau', X_i\leq\tau\}} \right|
= o_p(1).
$
Next, letting $f_{i-1}(\theta_1) = \frac{\dB b_{i-1}^* \dA \sigma_ 
{i-1} }{\sigma_{i-1}^{3}}$,
we estimate that
for $l>0$
\begin{eqnarray*}
& & E \left| \frac{1}{\sqrt{n h_n}} \sum_{i=1}^n
\sup_{\theta_1}
\left| \delta_{\theta_1} \delta_{\theta_2} \ell(i,i-1;\theta_1,  
\theta_{2,0}) \right| \right|^{2 l} \\
&\leq& \frac{C}{(n h_n)^l} E \left[ \sum_{i=1}^n \int_{t_{i-1}}^{t_i}  
\sup_{\theta_1} |f_{i-1}(\theta_1)|
\sigma(X_s,\theta_{1,0}) dW_s \right]^{2 l} \\
&& + \frac{C}{(n h_n)^l} E \left[ \sum_{i=1}^n \int_{t_{i-1}}^{t_i}  
\sup_{\theta_1} |f_{i-1}(\theta_1)|
(b(X_s,\theta_{2,0}) - b_{i-1}^* )ds \right]^{2 l} \\
&\leq& \frac{C}{(n h_n)^l} (n h_n)^{l-1} \sum_{i=1}^n  E \left[ \int_ 
{t_{i-1}}^{t_i} \sup_{\theta_1} |f_{i-1}(\theta_1)|^{2l}
\sigma^{2 l} (X_s,\theta_{1,0}) ds \right] \\
&& + \frac{C}{(n h_n)^l} (n h_n)^{2l-1} \sum_{i=1}^n E \left[  \int_ 
{t_{i-1}}^{t_i} \sup_{\theta_1} |f_{i-1}(\theta_1)|^{2 l}
(b(X_s,\theta_{2,0}) - b_{i-1}^* )^{2 l} ds \right] \\
&=& O(1).
\end{eqnarray*}
Consequently, one has that $|B_n| = o_p(1)$.



Next, we will prove that
\begin{equation} \label{AN}
\bar{\cal L}_n
\overset{d}{\to}
N(0, 4\Sigma).
\end{equation}
Let
\begin{eqnarray*}
\xi_i^{(1)}
&=& \frac{1}{\sqrt{n}}\dA \ell(i,i-1;\theta_{1,0})\chi_{\{X_{i-1}> \tau'\}} \\
&=& \frac{2}{\sqrt{n}}  \frac{\dA\sigma_{i-1}^*}{\sigma_{i-1}^*} \left(
1 - \frac{\left(\Delta_i X\right)^2}{h_n{{\sigma_{i-1}^{*2}}}}
\right)\chi_{\{X_{i-1}>\tau'\}},\\
\xi_i^{(2)} &=&  \frac{1}{\sqrt{nh_n}}\dB \ell(i,i-1;\theta_0)\chi_{\{X_{i-1}>\tau'\}} \\
&=&  -\frac{2}{\sqrt{n h_n}} \left\{\dB b_{i-1}^* \frac{\Delta_i X -  
b_{i-1}^* h_n}{\sigma_{i-1}^{*2}}
\right\}\chi_{\{X_{i-1}>\tau'\}}, \\
I(\theta_0) &=&
\begin{pmatrix}
I^{(1,1)}(\theta_{1,0}) & 0 \\
0 & I^{(2,2)}(\theta_0)
\end{pmatrix}
:= 4 \Sigma.
\end{eqnarray*}
In order to obtain (\ref{AN}),
by the combination of Theorems 3.2 and 3.4 of Hall and {Heyde} (1980),
it is enough to prove the following convergences.

\begin{eqnarray}
&\sum\limits_{i=1}^n& \E_{\theta_0}\left\{ \xi_{i}^{(k)} | \mathcal F_ 
{i-1}\right\} \overset{p}{\to} 0, \quad k=1,2, \label{eq:xi1}\\
&\sum\limits_{i=1}^n& \E_{\theta_0}\left\{ \left(\xi_{i}^{(k)} \right) 
^2 | \mathcal F_{i-1}\right\} \overset{p}{\to} I^{(k,k)}, \quad  
k=1,2,  \label{eq:xi2}\\
&\sum\limits_{i=1}^n& \E_{\theta_0}\left\{ \xi_{i}^{(1)} \xi_{i}^ 
{(2)} | \mathcal F_{i-1}\right\} \overset{p}{\to} 0, \label{eq:xi3}\\
&\sum\limits_{i=1}^n& \E_{\theta_0}\left\{ \left(\xi_{i}^{(k)} \right) 
^4 | \mathcal F_{i-1}\right\} \overset{p}{\to} 0, \quad k=1,2.  \label 
{eq:xi4}
\end{eqnarray}
For the proof of (\ref{eq:xi1}),
by using the It\^o-Taylor expansion and (\ref{I-T1}), one has
$$
\begin{aligned}
\sum\limits_{i=1}^n \E_{\theta_0}\left\{ \xi_{i}^{(1)} | \mathcal F_{i-1}\right\}
&=  \sqrt{n h_n^2} \cdot \frac{1}{n} \sum\limits_{i=1}^n  R(1,X_ 
{i-1}) \overset{p}{\to}  0.
\end{aligned}
$$
Moreover, since
$$
\E_{\theta_0}(X_i - X_{i-1}  |\mathcal F_{i-1}) = h_n b_{i-1}^* + R 
(h_n^2, X_{i-1}),
$$
we have
$$
\begin{aligned}
\sum_{i=1}^n \E_{\theta_0}( \xi_{i}^{(2)}  |\mathcal F_{i-1})\chi_{\{X_{i-1}>\tau'\}}
&=\frac{-2 \sqrt{n h_n^3} }{n}\sum_{i=1}^n R(1, X_{i-1}) \overset{p} 
{\to}  0,
\end{aligned}
$$
which completes the proof of (\ref{eq:xi1}).
For the proof of (\ref{eq:xi2}), noting that
\begin{eqnarray*}
& & E \left\{ \left(1 - \frac{\left(\Delta_i X \right)^2}{\sigma_{i-1} 
^{*2} h_n}
\right)^2 | \mathcal F_{i-1}\right\}  \\
&=&
1 + \frac{3 h_n^2 \sigma_{i-1}^{*4} + R(h_n^{5/2},X_{i-1})}{\sigma_ 
{i-1}^{*4} h_n^2} -2 \frac{h_n \sigma_{i-1}^{*2} + R(h_n^2, X_{i-1})} 
{\sigma_{i-1}^{*2} h_n} \\
&=& 2 + \sqrt{h_n} R(1, X_{i-1}),
\end{eqnarray*}
one has
$$
\begin{aligned}
\sum_{i=1}^n \E_{\theta_0}\left\{ \left(\xi_{i}^{(1)} \right)^2 |  
\mathcal F_{i-1}\right\}
&=
\sum_{i=1}^n \frac{4}{n}  \left(\frac{\dA\sigma_{i-1}^*}{\sigma_{i-1} 
^*}\right)^2
(2 + \sqrt{h_n} R(1, X_{i-1}))
  \chi_{\{X_{i-1}>\tau'\}} \\
&  \overset{p}{\to}  I^{(1,1)}(\theta_{1,0}),
\end{aligned}
$$
which proves \eqref{eq:xi2} for $k=1$.
It follows from the It\^o-Taylor expansion of
$\E_{\theta_0}\{( X_{i} - X_{i-1} - h_n b_{i-1}^*)^2| \mathcal F_{i-1} 
\}$ that
$$
\begin{aligned}
\sum_{i=1}^n \E_{\theta_0}&\left\{\left( \xi_i^{(2)} \right)^2|  
\mathcal F_{i-1}\right\} \\
&=\frac{1}{n h_n} 4 \sum_{i=1}^n \frac{(\dB b_{i-1}^*)^2}{\sigma^{*4}_ 
{i-1}} (h_n \sigma_{i-1}^{*2} + R(h_n^2, X_{i-1}))\chi_{\{X_{i-1}> 
\tau'\}}
  \\
&  \overset{p}{\to}  I^{(2,2)}(\theta_{0})
\end{aligned}
$$
and \eqref{eq:xi2} is proved.
For the proof of (\ref{eq:xi3}),
we consider
$$
\begin{aligned}
  \xi_i^{(1)} \xi_i^{(2)}
=&
-\frac{4}{n \sqrt{h_n}}  \frac{\dA\sigma_{i-1}^* \dB b_{i-1}^*} 
{\sigma_{i-1}^{*3}} \\
&\times\left(1 - \frac{\left(\Delta_i X  \right)^2}{\sigma_{i-1}^{*2}  
h_n}
\right)  \left\{ \Delta_i X - b_{i-1}^* h_n
\right\}\chi_{\{X_{i-1}>\tau'\}}.
\end{aligned}
$$
Since
$$
  \E_{\theta_0}\left\{\left(\Delta_i X \right)^2\left(\Delta_i X - b_ 
{i-1}^* h_n\right) | \mathcal F_{i-1}\right\}
  = R(h_n^2, X_{i-1})
$$
and
$$
  \E_{\theta_0}\left\{ \Delta_i X - b_{i-1}^* h_n| \mathcal F_{i-1} 
\right\}
  = R(h_n^2, X_{i-1}),
$$
one has
$$
\begin{aligned}
\sum_{i=1}^n
  \E_{\theta_0}\left\{  \xi_i^{(1)} \xi_i^{(2)} | \mathcal F_{i-1} 
\right\}
=&
-\frac{4}{n}  \sum_{i=1}^n \frac{\dA\sigma_{i-1}^* \dB b_{i-1}^*} 
{\sigma_{i-1}^{*3}}\chi_{\{X_{i-1}>\tau'\}}\\
&\times
\frac{1}{\sqrt{h_n}}\left( h_n^2 R(1, X_{i-1}) - \frac{h_n R(1, X_ 
{i-1})}{\sigma_{i-1}^{*2}}\right) \\
\overset{p}{\to}& \ 0.
\end{aligned}
$$
Hence \eqref{eq:xi3} is proved.
For the proof of (\ref{eq:xi4}),
using the estimate
$$\E_{\theta_0}\{\left(\Delta_i X \right)^{2k}| \mathcal F_{i-1}\} =  
h_n^{k} R(1,X_{i-1}),$$
one has
$$
\E_{\theta_0}\left\{ \sum_{i=1}^n \left(\xi_i^{(1)}\right)^4 |  
\mathcal F_{i-1}\right\} \leq \frac{C'}{n}
\frac{1}{n} \sum_{i=1}^n \left(\frac{\dA\sigma_{i-1}^*}{\sigma_{i-1} 
^*}\right)^4
  \chi_{\{X_{i-1}>\tau'\}} \{ 1+ R(1,X_{i-1}) \}
\overset{p}{\to}  0,
$$
which
completes the proof of \eqref{eq:xi4} for $k=1$.
For the case $k=2$,
by using the following estimate
$$
  \E_{\theta_0}\left\{  \left(\Delta_i X - b_{i-1}^* h_n\right)^4 |  
\mathcal F_{i-1}\right\} =  h_n^2 R(1,X_{i-1}),
  $$
we have that
$$
\sum_{i=1}^n  \E_{\theta_0}\left\{\left(\xi_i^{(2)} \right)^4 |  
\mathcal F_{i-1}\right\}
\leq \frac{C'}{n}\frac{1}{n} \sum_{i=1}^n\left(\frac{\dB b_{i-1}^*} 
{\sigma_{i-1}^{*2}}\right)^4
  \chi_{\{X_{i-1}>\tau'\}} R(1,X_{i-1})
\overset{p}{\to}  0.
$$
Thus \eqref{eq:xi4} is proved. This completes the proof of (\ref{AN}).
It follows from (\ref{score1}), (\ref{score2}) and (\ref{AN}) that
\begin{equation} \label{AN1}
{\cal L}_n
\overset{d}{\to}
N(0, 4\Sigma).
\end{equation}



Next we consider asymptotic properties of the observed information.
Let
$$
D_n(\theta) =
\left(
\begin{array}{cc}
\frac{1}{{n}} \delta_{\theta_1}^2 g_n(\theta_{1}) & 0 \\
(\theta)
0
& \frac{1}{n h_n} \delta_{\theta_2}^2 \ell_n(\hat{\theta}_{1,n},  
\theta_2)
\end{array}
\right), \quad
D(\theta) =
\left(
\begin{array}{cc}
\bar{\cal G}(\theta_1)  & 0 \\
0 & \bar{\cal L}(\theta_2)
\end{array}
\right),
$$
where
\begin{eqnarray*}
\bar{\cal G}(\theta_1) &=&
2 \int_{\bf R} \frac{\delta_{\theta_1}^2 \sigma(X,\theta_1)}{\sigma^3 
(x,\theta_1)} \left( \sigma^2(x,\theta_1) -\sigma^2(x,\theta_{1,0})  
\right)
\chi_{\{ x >\tau \}} \nu_{\theta_0}(dx) \\
& &+ 2 \int_{\bf R} \frac{\left( 3 \sigma^2(x,\theta_{1,0}) - \sigma^2 
(x,\theta_1) \right)
\left( \delta_{\theta_1} \sigma(x,\theta_1) \right)^2}{\sigma^4(x, 
\theta_1)} \chi_{\{ x >\tau \}} \nu_{\theta_0}(dx), \\
\bar{\cal L}(\theta_2) &=&
2 \int_{\bf R} \left( \frac{\dB b(x,\theta_2)}{\sigma(x,\theta_ 
{1,0})} \right)^2 \chi_{\{ x > \tau \}} \nu_{\theta_0}(dx)  \\
& & - 2  \int_{\bf R}
\frac{( b(x,\theta_{2,0}) - b(x,\theta_2 ) ) \dBB b(x,\theta_2) } 
{\sigma^2(x,\theta_{1,0})}
\chi_{\{ x > \tau \}} \nu_{\theta_0}(dx).
\end{eqnarray*}
In order to prove that
\begin{eqnarray} \label{OI}
\sup_\theta \left| D_n(\theta) - D(\theta) \right| &=& o_p(1), 
\end{eqnarray}
it is sufficient to show that
\begin{eqnarray}
\sup_{\theta_1} \left| \frac{1}{n} \delta_{\theta_1}^2 g_n(\theta_1) -
\frac{1}{n} \delta_{\theta_1}^2 \bar{g}_n(\theta_1)  \right|
&=& o_p(1), \label{uchida2}
\\
\sup_\theta \left|
\frac{1}{{n h_n}} \delta_{\theta_2}^2 \ell_n(\theta) -
\frac{1}{{n h_n}} \delta_{\theta_2}^2 \bar\ell_n(\theta)
\right|
&=& o_p(1), \label{uchida3}
\\
\sup_{\theta_1} \left| \frac{1}{n} \delta_{\theta_1}^2 \bar{g}_n 
(\theta_1) -
\bar{\cal G}(\theta_1)  \right|
&=& o_p(1), \label{uchida4}
\\
\sup_{\theta_2} \left|
\frac{1}{{n h_n}} \delta_{\theta_2}^2 \bar\ell_n(\hat{\theta}_{1,n},  
\theta_2) -
\bar{\cal L}(\theta_2)
\right|
&=& o_p(1). \label{uchida5}
\end{eqnarray}
For the proof of (\ref{uchida2}),
as in the proof of the uniform convergence of (\ref{uchida:cont1-b}),
one has that 
\begin{eqnarray*}
& & \E \left\{ \sup_{\theta_1} \left|\frac{1}{n} \sum_{i=1}^n  \delta_ 
{\theta_1}^2 g(i, i-1; \theta_1)  \chi_{\{X_{i-1}>\tau', X_i\leq\tau 
\}} \right| \right\} \\
&\leq&
\frac{1}{n} \sum_{i=1}^n \left|\left| \ \sup_{\theta_1} \left| \delta_ 
{\theta_1}^2 g(i, i-1; \theta_1) \right|\,\,
\right|\right|_2
  P(X_{i-1}>\tau', X_i\leq\tau)^\frac{1}{2} \\
&\rightarrow& 0.
\end{eqnarray*}
For the proof of (\ref{uchida3}),
in a quite similar way as in the proof of (\ref{uchida2}),
one has that 
\begin{eqnarray*}
& & \E \left\{ \sup_{\theta} \left|\frac{1}{n h_n } \sum_{i=1}^n   
\delta_{\theta_1}^2
\ell(i, i-1; \theta)  \chi_{\{X_{i-1}>\tau', X_i\leq\tau\}} \right|  
\right\} \\
&\leq&
\frac{1}{n h_n } \sum_{i=1}^n \left|\left| \ \sup_{\theta} \left|  
\delta_{\theta_1}^2 \ell(i, i-1; \theta) \right|\,\,
\right|\right|_2
  P(X_{i-1}>\tau', X_i\leq\tau)^\frac{1}{2} \\
&\leq& C \frac{1}{ h_n^{1/2}} \times \left( \frac{h_n^\alpha}{\tau'- 
\tau} \right)^{k/2} (h_n^{1/2-\alpha})^{k/2} \\
&=& O \left( h_n^{k/4- \alpha k /2 -1/2} \right) \rightarrow 0,
\end{eqnarray*}
where we took $k > 2/(1- 2 \alpha)$ in (\ref{yoshida:ess2}).
For the proof of (\ref{uchida4}),
we set
$$
\begin{aligned}
\eta_i(\theta_1) =& \frac{1}{n} \dAA \bar{g}(i,i-1;\theta_{1})\chi_{\{X_{i-1}>\tau'\}}\\
=& \frac{2}{n h_n \sigma_{i-1}^{4}}  \biggl\{(3 (\Delta_i X)^2 - h_n  
\sigma_{i-1}^{2}) (\dA \sigma_{i-1})^{2} \\
&+ \sigma_{i-1}(h_n \sigma_{i-1}^{2} -  (\Delta_i X)^2 )\dAA \sigma_ 
{i-1}\biggr\}\chi_{\{X_{i-1}>\tau'\}}.
\end{aligned}
$$
It follows from standard arguments that
$$
\begin{aligned}
\sum_{i=1}^n \E_{\theta_0}\left\{ \eta_i(\theta_1) | \mathcal F_{i-1} 
\right\}
&\overset{p}{\to}
\bar{\cal G}(\theta_1), \\
\sum_{i=1}^n \E_{\theta_0}\left\{ (\eta_i(\theta_1))^2 | \mathcal F_ 
{i-1}\right\}
&\overset{p}{\to} 0.
\end{aligned}
$$
Therefore one has that for each $\theta_1$,
$$\frac{1}{n} \dAA \bar{g}_n(\theta_1) \overset{p}{\to} \bar{\cal G} 
(\theta_1).$$
It is easy to show that $\sup_n E[ \sup_{\theta_1} | \frac{1}{n} 
\delta_{\theta_1}^3  \bar{g}_n(\theta_1) | ] < \infty$,
which completes the proof of (\ref{uchida4}).
For the proof of (\ref{uchida5}), we set
\begin{eqnarray*}
\frac{1}{n h_n} \dBB \bar{\ell}_n(i,i-1; \hat{\theta}_{1,n},  
\theta_2) 
&=& \Xi_1(\theta_2) + \Xi_2(\theta_2) + \Xi_3(\theta_2),
\end{eqnarray*}
where
\begin{eqnarray*}
\Xi_1(\theta_2)
&=& \frac{2}{n} \sum_{i=1}^n \left\{ \left(\frac{\dB b_{i-1}}{\hat 
{\sigma}_{i-1}}\right)^2 -
\frac{(b_{i-1}^* - b_{i-1} )\dBB b_{i-1} }{\hat{\sigma}_{i-1}^{2}}  
\right\} \chi_{\{X_{i-1}>\tau'\}}, \\
\Xi_2(\theta_2)
&=&
- \frac{2}{n h_n} \sum_{i=1}^n
\frac{\dBB b_{i-1} \int_{t_{i-1}}^{t_i} \{ b(X_s,\theta_{2,0})-b_{i-1} 
^* \} ds }{\hat{\sigma}_{i-1}^{2}}
\chi_{\{X_{i-1}>\tau'\}}, \\
\Xi_3(\theta_2)
&=&
- \frac{2}{n h_n} \sum_{i=1}^n
\frac{\dBB b_{i-1} \int_{t_{i-1}}^{t_i} \sigma(X_s,\theta_{1,0})  
dW_s }{\hat{\sigma}_{i-1}^{2}}
\chi_{\{X_{i-1}>\tau'\}}.
\end{eqnarray*}
In a quite similar way as in the proof of (\ref{convT1}),
\begin{eqnarray*}
\sup_{\theta_2} \left| \Xi_1(\theta_2) - \bar{\cal L}(\theta_2)  
\right| &=& o_p(1), \\
\sup_{\theta_2} \left| \Xi_2(\theta_2)  \right| &=& o_p(1), \\
\sup_{\theta_2} \left| \Xi_3(\theta_2)  \right| &=& o_p(1).
\end{eqnarray*}
This completes the proof of (\ref{uchida5}).
Thus, (\ref{OI}) is proved.

By the Taylor expansion,
$$
\int_0^1 D_n(\theta_0 + u(\hat{\theta}_n - \theta_0)) du S_n = {\cal  
L}_n
$$
on an event with probability tending to one,
where
$$
S_n =
\left(
\begin{array}{c}
{\sqrt{n}} (\hat{\theta}_{1,n} - \theta_{1,0}) \\
{\sqrt{n h_n}} (\hat{\theta}_{2,n} - \theta_{2,0})
\end{array}
\right).
$$
It follows from (\ref{AN1}) that
\begin{equation} \label{uchida-asy1}
{\cal L}_n \overset{d}{\to} N(0, 4 \Sigma).
\end{equation}
By (\ref{OI}) and the continuity of $D(\theta)$ with respect to $ 
\theta$,
one has
\begin{eqnarray}
D_n(\theta_0)  &\overset{p}{\to}& 2 \Sigma, \label{uchida-asy2} \\
\sup_{|\theta| \leq \epsilon_n} \left| D_n(\theta_0 + \theta) - D_n 
(\theta_0) \right| &=& o_p(1) \label{uchida-asy3}
\end{eqnarray}
{{for any sequence $\epsilon_n$ of positive numbers
tending to zero.}}
By using (\ref{uchida-asy1})-(\ref{uchida-asy3}), it is easy to  
obtain the desired result.
This completes the proof.
\end{proof}

\section*{Acknowledgments}
This work was done during the period of staying of  the first author  
at the ISM (Institute of Statistical Mathematics)
and the University of Tokyo. The ISM is acknowledged with thanks.
The works have been supported by the JSPS (Japan Society for the  
Promotion of Science) Program FY2006, grant ID No. S06174.
The researches of Masayuki Uchida
and Nakahiro Yoshida were supported by
Grants-in-Aid for Scientific Research from the JSPS,
and by Cooperative Research Program
of the Institute of Statistical Mathematics.

\end{document}